\documentclass[11pt]{article}
\usepackage[latin1]{inputenc}
\usepackage[T1]{fontenc}
\usepackage[english]{babel}
\usepackage{amsmath}
\usepackage{amsthm}
\usepackage{amsfonts}
\usepackage{amssymb}
\usepackage{graphicx}
\usepackage{stmaryrd}
\usepackage{tikz}
\usepackage{url}
\usepackage[np]{numprint}
\usepackage{vmargin}
\usepackage{mathrsfs}
\newtheorem{thm}{Theorem}
\newtheorem{lem}{Lemma}
\newtheorem{prop}{Proposition}
\newtheorem{coro}{Corollary}
\theoremstyle{definition}
\newtheorem{defi}{Definition}
\theoremstyle{remark}
\newtheorem{rem}{Remark}

\newcommand{\abs}[1]{\left\lvert#1\right\rvert}

\setmarginsrb{2cm}{2cm}{2cm}{2.5cm}{0cm}{0cm}{0cm}{1cm}
\newcounter{ex}

\begin{document}

\date{}
\title{\textbf{The Calkin-Wilf tree of a quadratic surd}}
\author{Lionel Ponton \\
%EndAName
{\small \texttt{lionel.ponton@gmail.com}}}
\maketitle

\begin{abstract}
By using the Calkin-Wilf tree, we prove the irrationality of numbers of the form $\alpha=\frac{\sqrt{N}+p}{q}$ where $N$ is a positive integer which is not a perfect square, $p$ is a rational integer such that $p^2<N$ and $q$ is a positive integer which divides $N-p^2$. For this, we consider an analogue of the Calkin-Wilf tree with root $\alpha$ and we define a special path in this tree which satisfies remarkable properties of periodicity and symmetry. This path is closely related to the continued fraction expansion of $\alpha$ and allows us to give new proofs of theorems due to Legendre and to Galois about the form of such an expansion in special cases of square roots and reduced quadratic surds.
\end{abstract}

\bigskip

\bigskip

\section{Introduction}

In 2000, N. Calkin and H. S. Wilf \cite{CW00} have defined a binary tree in the following way:

\begin{enumerate}
\item[$\bullet $] the root of the tree is $\frac{1}{1}$;

\item[$\bullet $] the vertex labeled $\frac{a}{b}$ has two children: the left child labeled $\frac{a}{a+b}$ and the right child labeled $\frac{a+b}{b}$.
\end{enumerate}

In other words, the two children of a positive rational number $x$ are $\frac{x}{x+1}$ and $x+1$:

\begin{center}
\begin{tikzpicture}[xscale=1,yscale=1]
\node (P) at ({0},{0}) {$x$};
\node (fg) at ({-2},{-1}) {$\dfrac{x}{x+1}$};
\node (fd) at ({2},{-1}) {$x+1$};
\draw (P)--(fg);
\draw (P)--(fd);
\end{tikzpicture}
\end{center}

Thus, the first few rows of the Calkin-Wilf tree are:

\begin{center}
\begin{tikzpicture}[xscale=1,yscale=1]
\node (P) at ({0},{0}) {$1$};
\node (fg) at ({-4},{-1}) {$\dfrac{1}{2}$};
\node (fgg) at ({-6},{-2}) {$\dfrac{1}{3}$};
\node (fggg) at ({-7},{-3}) {$\dfrac{1}{4}$};
\node (fggd) at ({-5},{-3}) {$\dfrac{4}{3}$};
\node (fgd) at ({-2},{-2}) {$\dfrac{3}{2}$};
\node (fgdg) at ({-3},{-3}) {$\dfrac{3}{5}$};
\node (fgdd) at ({-1},{-3}) {$\dfrac{5}{2}$};
\node (fd) at ({4},{-1}) {$2$};
\node (fdg) at ({2},{-2}) {$\dfrac{2}{3}$};
\node (fdgg) at ({1},{-3}) {$\dfrac{2}{5}$};
\node (fdgd) at ({3},{-3}) {$\dfrac{5}{3}$};
\node (fdd) at ({6},{-2}) {$3$};
\node (fddg) at ({5},{-3}) {$\dfrac{3}{4}$};
\node (fddd) at ({7},{-3}) {$4$};
\draw (P)--(fg);
\draw (fg)--(fgg);
\draw (fgg)--(fggg);
\draw (fgg)--(fggd);
\draw (fg)--(fgd);
\draw (fgd)--(fgdg);
\draw (fgd)--(fgdd);
\draw (P)--(fd);
\draw (fd)--(fdg);
\draw (fdg)--(fdgg);
\draw (fdg)--(fdgd);
\draw (fd)--(fdd);
\draw (fdd)--(fddg);
\draw (fdd)--(fddd);
\end{tikzpicture}
\end{center}

In their article, Calkin and Wilf have shown that every positive rational number appears once and only once in reduced form in this tree. Therefore, by following the breadth-first order, one obtains an explicit enumeration of the positive rationals. The first few terms of this sequence, called the Calkin-Wilf sequence, are:

\begin{equation*}
1,~\dfrac{1}{2},~2,~\dfrac{1}{3},~\dfrac{3}{2},~\dfrac{2}{3},~3,~\dfrac{1}{4},~\dfrac{4}{3},~\dfrac{3}{5},~\dfrac{5}{2},~\dfrac{2}{5},~\dfrac{5}{3},~\dfrac{3}{4},~4,~\dfrac{1}{5},~\dfrac{5}{4},~\dfrac{4}{7},...
\end{equation*}

It has lots of remarkable properties and it is related to many others subjects such as Stern sequence, Farey sequences, continued fractions... (see \cite[Chap. 19]{AZ14} and \cite{Nor10}).

\bigskip

In this paper, we consider, for any real number $x$, the analogue of the Calkin-Wilf tree with root $x$. We call it the Calkin-Wilf tree of $x$. We are particularly interested in the case $x=\frac{\sqrt{N}+p}{q}$ where $N$ is a positive integer which is not a perfect square, $p$ is a rational integer such that $p^2<N$ and $q$ is a positive integer which divides $N-p^2$. In this case, the Calkin-Wilf tree of $x$ allows us to prove that $x$ is irrational by considering a special path which satisfies remarkable properties of periodicity and symmetry. A careful study of this path reveals close links with continued fraction expansion and leads us to give new proofs of theorems due to Legendre and to Galois. 

\section{Two examples}

\textbf{Irrationality of $\sqrt{2}$}

The first few levels of the Calkin-Wilf tree of $\sqrt{2}$ are:

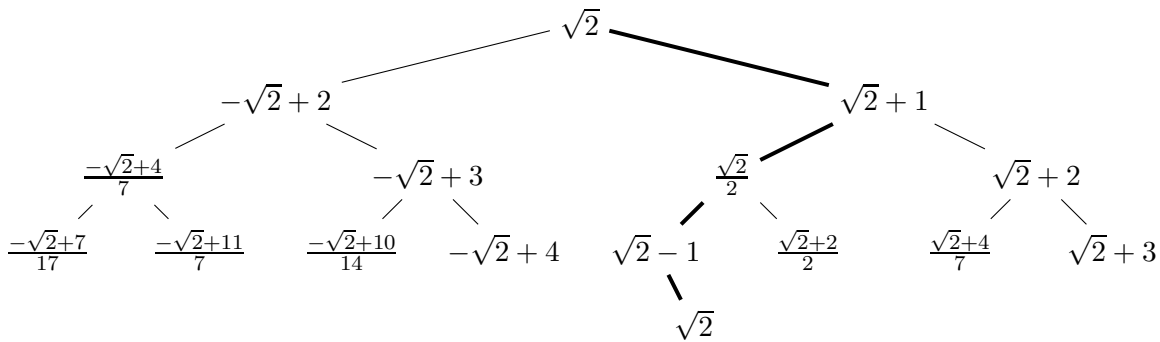
\begin{center}
\begin{tikzpicture}[xscale=1,yscale=1]
\node (P) at ({0},{0}) {$\sqrt{2}$};
\node (fg) at ({-4},{-1}) {$-\sqrt{2}+2$};
\node (fgg) at ({-6},{-2}) {$\frac{-\sqrt{2}+4}{7}$};
\node (fggg) at ({-7},{-3}) {$\frac{-\sqrt{2}+7}{17}$};
\node (fggd) at ({-5},{-3}) {$\frac{-\sqrt{2}+11}{7}$};
\node (fgd) at ({-2},{-2}) {$-\sqrt{2}+3$};
\node (fgdg) at ({-3},{-3}) {$\frac{-\sqrt{2}+10}{14}$};
\node (fgdd) at ({-1},{-3}) {$-\sqrt{2}+4$};
\node (fd) at ({4},{-1}) {$\sqrt{2}+1$};
\node (fdg) at ({2},{-2}) {$\frac{\sqrt{2}}{2}$};
\node (fdgg) at ({1},{-3}) {$\sqrt{2}-1$};
\node (fdgd) at ({3},{-3}) {$\frac{\sqrt{2}+2}{2}$};
\node (fdd) at ({6},{-2}) {$\sqrt{2}+2$};
\node (fddg) at ({5},{-3}) {$\frac{\sqrt{2}+4}{7}$};
\node (fddd) at ({7},{-3}) {$\sqrt{2}+3$};
\node (fdggd) at ({1.5},{-4}) {$\sqrt{2}$};
\draw (P)--(fg);
\draw (fg)--(fgg);
\draw (fgg)--(fggg);
\draw (fgg)--(fggd);
\draw (fg)--(fgd);
\draw (fgd)--(fgdg);
\draw (fgd)--(fgdd);
\draw[line width = 1.5pt] (P)--(fd);
\draw[line width = 1.5pt] (fd)--(fdg);
\draw[line width = 1.5pt] (fdg)--(fdgg);
\draw (fdg)--(fdgd);
\draw (fd)--(fdd);
\draw (fdd)--(fddg);
\draw (fdd)--(fddd);
\draw[line width = 1.5pt] (fdgg)--(fdggd);
\end{tikzpicture}
\end{center}

Assumes that $\sqrt{2}$ is rational. Then $\sqrt{2}$ appears somewhere in the Calkin-Wilf tree and so the Calkin-Wilf tree of $\sqrt{2}$ is a subtree of the Calkin-Wilf tree. But, we see that $\sqrt{2}$ appears twice in the Calkin-Wilf tree of $\sqrt{2}$: once at the top-level and once as the right child of $\sqrt{2}-1$ at the fourth level. This is impossible since the labels in the Calkin-Wilf tree are all different. We conclude that $\sqrt{2}$ is irrational.

\bigskip

One can remark two interesting facts:
\begin{enumerate}
\item[$\bullet$] it is easy to prove by induction that a positive rational $x$ and its reciprocal $\frac{1}{x}$ are at the same level in the Calkin-Wilf tree (see \cite[Theorem 1]{BBT10}). Since $\frac{\sqrt{2}}{2}$ is two levels below $\sqrt{2}$ in the Calkin-Wilf tree of $\sqrt{2}$, we obtain an other contradiction which yields the irrationality of $\sqrt{2}$.
\item[$\bullet$] In some proofs of the irrationality of $\sqrt{2}$, one shows, more or less explicitly, that if $\sqrt{2}=\frac{a}{b}$ is a rational in reduced form then $\frac{2b-a}{a-b}=\sqrt{2}$ which contradicts the minimality of $b$ since $0 < a-b < b$ (see \cite[Proof 8]{Bog}). By following the path from $\sqrt{2}$ to $\frac{\sqrt{2}}{2}$ in the above tree, it yields:
\[\sqrt{2}=\dfrac{a}{b} \longrightarrow \sqrt{2}+1=\dfrac{a+b}{b} \longrightarrow \dfrac{\sqrt{2}}{2}=\dfrac{a+b}{a+2b}\]
or, in the reverse order,
\[\dfrac{\sqrt{2}}{2}=\dfrac{b}{a} \longrightarrow \sqrt{2}+1=\dfrac{b}{a-b} \longrightarrow \sqrt{2}=\dfrac{2b-a}{a-b}\]
\end{enumerate}
and thus $\frac{2b-a}{a-b}=\frac{a}{b}$, which yields the same contradiction. 

\bigskip

\textbf{Irrationality of $\frac{\sqrt{3}+1}{2}$}

In the same way, the first three levels of the Calkin-Wilf tree of $\frac{\sqrt{3}+1}{2}$ are:

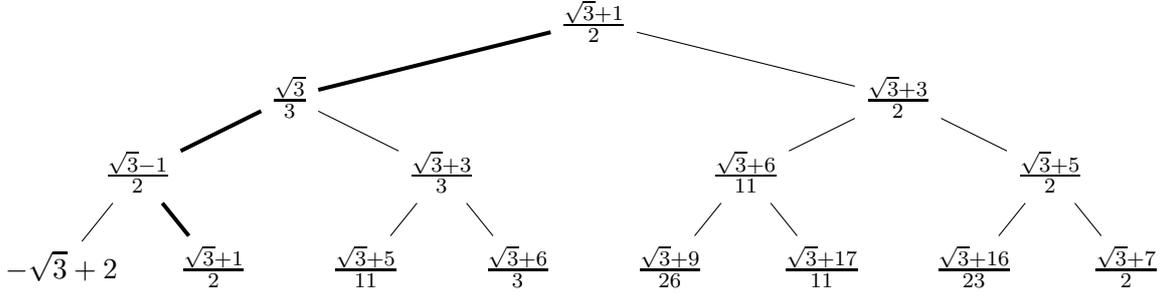
\begin{center}
\begin{tikzpicture}[xscale=1,yscale=1]
\node (P) at ({0},{0}) {$\frac{\sqrt{3}+1}{2}$};
\node (fg) at ({-4},{-1}) {$\frac{\sqrt{3}}{3}$};
\node (fgg) at ({-6},{-2}) {$\frac{\sqrt{3}-1}{2}$};
\node (fggg) at ({-7},{-3.25}) {$-\sqrt{3}+2$};
\node (fggd) at ({-5},{-3.25}) {$\frac{\sqrt{3}+1}{2}$};
\node (fgd) at ({-2},{-2}) {$\frac{\sqrt{3}+3}{3}$};
\node (fgdg) at ({-3},{-3.25}) {$\frac{\sqrt{3}+5}{11}$};
\node (fgdd) at ({-1},{-3.25}) {$\frac{\sqrt{3}+6}{3}$};
\node (fd) at ({4},{-1}) {$\frac{\sqrt{3}+3}{2}$};
\node (fdg) at ({2},{-2}) {$\frac{\sqrt{3}+6}{11}$};
\node (fdgg) at ({1},{-3.25}) {$\frac{\sqrt{3}+9}{26}$};
\node (fdgd) at ({3},{-3.25}) {$\frac{\sqrt{3}+17}{11}$};
\node (fdd) at ({6},{-2}) {$\frac{\sqrt{3}+5}{2}$};
\node (fddg) at ({5},{-3.25}) {$\frac{\sqrt{3}+16}{23}$};
\node (fddd) at ({7},{-3.25}) {$\frac{\sqrt{3}+7}{2}$};
\draw[line width = 1.5pt] (P)--(fg);
\draw[line width = 1.5pt] (fg)--(fgg);
\draw (fgg)--(fggg);
\draw[line width = 1.5pt] (fgg)--(fggd);
\draw (fg)--(fgd);
\draw (fgd)--(fgdg);
\draw (fgd)--(fgdd);
\draw (P)--(fd);
\draw (fd)--(fdg);
\draw (fdg)--(fdgg);
\draw (fdg)--(fdgd);
\draw (fd)--(fdd);
\draw (fdd)--(fddg);
\draw (fdd)--(fddd);
\end{tikzpicture}
\end{center}

Since $\frac{\sqrt{3}+1}{2}$ appears twice, we conclude similarly that $\frac{\sqrt{3}+1}{2}$ is irrational.

\section{General case}

Throughout this paper, $N$ is a positive integer which is not a perfect square, $p$ is a rational integer such that $p^2<N$ and $q$ is a positive integer which divides $N-p^2$. Furthermore, we put
\[\alpha:=\dfrac{\sqrt{N}+p}{q}.\]
Note that $\alpha>0$ and so, by definition, all the numbers in the Calkin-Wilf tree of $\alpha$ are positive too.

As in previous proofs, to derive the irrationality of $\alpha$, it suffices to show that the Calkin-Wilf tree of $\alpha$ contains the same number twice. We can note that in the above examples the repeated number lies in a path descending from the root of the tree. We call this kind of path an \emph{irrationality path}. We will now show that one can always find an irrationality path -- which proves the irrationality of $\alpha$ -- and we will describe a very simple algorithm to find it.

\begin{prop} --- Every number in the Calkin-Wilf tree of $\alpha$ can be written in the form 
\[\dfrac{\varepsilon\sqrt{N}+c}{d}\]
where $\varepsilon\in\{-1, 1\}$, $c$ is a rational integer and $d$ is a positive integer which divides $N-c^2$.
\label{prop_e_c_d}
\end{prop}

\emph{Proof}. --- The property is true at the top-level with $\varepsilon=1$, $c=p$ and $d=q$.

Assume that any number $x$ at a certain level $n$ in the tree can be written 
\[x=\dfrac{\varepsilon\sqrt{N}+c}{d}\]
where $\varepsilon\in\{-1, 1\}$, $c$ is a rational integer and $d$ is a positive integer such that $d \mid N-c^2$. 

The right child of $x$ is
\[x_r=x+1=\dfrac{\varepsilon\sqrt{N}+c+d}{d}=\dfrac{\varepsilon_r\sqrt{N}+c_r}{d_r}\]
where $\varepsilon_r:=\varepsilon\in\{-1,1\}$, $c_r:=c+d$ is a rational integer and $d_r:=d$ is a positive integer. Furthermore, $N-c_r^2=N-c^2-d(2c+d)$ is divisible by $d_r=d$ since $d$ divides $N-c^2$.

The left child of $x$ is 
\[x_{\ell}=\dfrac{x}{x+1}=\dfrac{\varepsilon\sqrt{N}+c}{\varepsilon\sqrt{N}+c+d}=\dfrac{\left(\varepsilon\sqrt{N}+c\right)\left(\varepsilon\sqrt{N}-c-d\right)}{N-(c+d)^2}=\dfrac{-\varepsilon d\sqrt{N}+N-c^2-cd}{N-(c+d)^2}.\]
with $N-(c+d)^2\neq 0$ since $N$ is not a perfect square.
If $N-(c+d)^2>0$, we put
\[\varepsilon_{\ell}:=-\varepsilon, \qquad c_{\ell}:=\dfrac{N-c^2}{d}-c \qquad \text{and} \qquad d_{\ell}:=\dfrac{N-(c+d)^2}{d}=\dfrac{N-c^2}{d}-2c-d\]
and otherwise we put
\[\varepsilon_{\ell}:=\varepsilon, \qquad c_{\ell}:=c-\dfrac{N-c^2}{d} \qquad \text{and} \qquad d_{\ell}:=\dfrac{(c+d)^2-N}{d}=2c+d-\dfrac{N-c^2}{d}.\]
In both cases, $x_{\ell}=\frac{\varepsilon_{\ell}\sqrt{N}+c_{\ell}}{d_{\ell}}$ where $\varepsilon_{\ell}\in\{-1,1\}$, $c_{\ell}$ is a rational integer and $d_{\ell}$ is a positive integer (since $d$ divides $N-c^2$). Moreover, 
\[N-c_{\ell}^2=N-c^2+2c\dfrac{N-c^2}{d}-\left(\dfrac{N-c^2}{d}\right)^2=\dfrac{N-c^2}{d}\left(d+2c-\dfrac{N-c^2}{d}\right)=\dfrac{N-c^2}{d} (\pm d_{\ell})\]
and, since $d$ divides $N-c^2$, $d_{\ell}$ divides $N-c_{\ell}^2$.

Hence, the property is true at level $n+1$ and Proposition \ref{prop_e_c_d} is proved by induction. \hfill $\square$

\begin{rem} \hspace{1cm}
\begin{enumerate}
\item We do not claim yet that the previous form is unique (even if it is true as we will see in Lemma \ref{lem_uniqueness}) so when we refer to $\varepsilon$, $c$ or $d$, one should understand the values as they are defined in the constructive proof of Proposition \ref{prop_e_c_d}.
\item One can notice the last computation in the previous proof shows that $N-c_{\ell}^2=\frac{N-c^2}{d}d_{\ell}$ if $\varepsilon_{\ell}=\varepsilon$.
\end{enumerate}
\label{rem_e_c_d}
\end{rem}

\begin{defi} --- We define a path in the Calkin-Wilf tree of $\alpha$ by the following algorithm.

\begin{center}
\begin{tabular}{l}
$\blacktriangleright$ Start the path at the top-level. \\
$\blacktriangleright$ From a vertex labeled $x$, \\
\begin{tabular}{l|l}
\phantom{mm}& If the value of $\varepsilon$ for the left child of $x$ is positive \\
&\hspace{0.5cm} go to the left \\
& Else \\
& \hspace{0.5cm} go to the right. \\
&End if
\end{tabular}
\end{tabular}
\end{center}

We call it the \emph{left-positive path of $\alpha$} (abbreviated LPP of $\alpha$).
\end{defi}

\begin{rem} \hspace{1cm}

\begin{enumerate}
\item Since a step right keeps the value of $\varepsilon$, $\varepsilon=1$ for every number in the LPP of $\alpha$. 
\item If $p=0$ and $N>q^2$ then $N-(p+q)^2>0$ so the LPP of $\alpha$ begins with a step to the right.  
\end{enumerate}
\label{rem_positive_path}
\end{rem}

For example, the beginning of the LPP of $\sqrt{5}$ is:

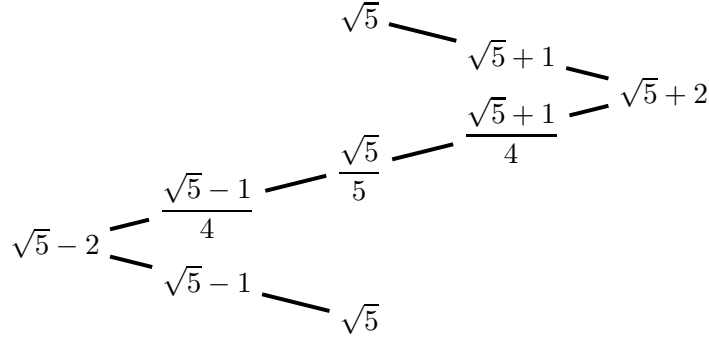
\begin{center}
\begin{tikzpicture}[xscale=1,yscale=1]
\node (n1) at ({0},{0}) {$\sqrt{5}$};
\node (n2) at ({2},{-0.5}) {$\sqrt{5}+1$};
\node (n3) at ({4},{-1}) {$\sqrt{5}+2$};
\node (n4) at ({2},{-1.5}) {$\dfrac{\sqrt{5}+1}{4}$};
\node (n5) at ({0},{-2}) {$\dfrac{\sqrt{5}}{5}$};
\node (n6) at ({-2},{-2.5}) {$\dfrac{\sqrt{5}-1}{4}$};
\node (n7) at ({-4},{-3}) {$\sqrt{5}-2$};
\node (n8) at ({-2},{-3.5}) {$\sqrt{5}-1$};
\node (n9) at ({0},{-4}) {$\sqrt{5}$};
\draw[line width=1.5pt] (n1)--(n2);
\draw[line width=1.5pt] (n2)--(n3);
\draw[line width=1.5pt] (n3)--(n4);
\draw[line width=1.5pt] (n4)--(n5);
\draw[line width=1.5pt] (n5)--(n6);
\draw[line width=1.5pt] (n6)--(n7);
\draw[line width=1.5pt] (n7)--(n8);
\draw[line width=1.5pt] (n8)--(n9);
\end{tikzpicture}
\end{center}

\bigskip

We can note that this path is an irrationality path for $\sqrt{5}$.

\bigskip

In the following, we denote $(x_n)_{n\geqslant 0}$ the sequence of numbers in the LPP of $\alpha$ in order of appearance and $(c_n)_{n\geqslant 0}$ and $(d_n)_{n\geqslant 0}$ the sequences of integers defined by Proposition \ref{prop_e_c_d} such that 
\[x_n=\frac{\sqrt{N}+c_n}{d_n}\]
for every non negative integer $n$. Note that, by definition, for every non negative integer $n$,
\begin{equation}
x_{n+1}=\dfrac{\sqrt{N}+c_n+d_n}{d_n} \text{ with } N-(c_n+d_n)^2>0 \text{ if $x_{n+1}$ is the right child of $x_n$}
\label{eq_x_n+1_right}
\end{equation}
and
\begin{equation}
x_{n+1}=\dfrac{\sqrt{N}+c_n-\frac{N-c_n^2}{d_n}}{\frac{(c_n+d_n)^2-N}{d_n}} \text{ with } (c_n+d_n)^2-N>0 \text{ if $x_{n+1}$ is the left child of $x_n$}.
\label{eq_x_n+1_left}
\end{equation}

Furthermore, for every non negative integer $n$, we put $s_n:=\ell$ if $x_{n+1}$ is the left child of $x_{n}$ and $s_n:=r$ if $x_{n+1}$ is the right child of $x_{n}$. Thus, the sequence $(s_n)$ encodes the successive steps in the LPP of $\alpha$. For example, if $\alpha=\sqrt{5}$, one has

\[x_0=\sqrt{5},~x_1=\sqrt{5}+1,~x_2=\sqrt{5}+2,~x_3=\dfrac{\sqrt{5}+1}{4},~x_4=\dfrac{\sqrt{5}}{5},~x_5=\dfrac{\sqrt{5}-1}{4}, ...\]
The first few terms of $(c_n)$ are
\[0,~1,~2,~1,~0,~-1,~-2,~-1,~0, ...\]
the first few terms of $(d_n)$ are 
\[1,~1,~1,~4,~5,~4,~1,~1,~1,  ...\]
and the first few terms of $(s_n)$ are
\[r,~r,~\ell,~\ell,~\ell,~\ell,~r,~r, ...\] 

\begin{prop} --- For every non negative integer $n$, $N>c_n^2$.
\label{prop_N_c_n}
\end{prop}

\emph{Proof}. --- Since $c_0=p$, the property is true for $n=0$ by definition.

Assume that the property is true for a certain non negative integer $n$.  

If $x_{n+1}$ is the right child of $x_n$ then, by \eqref{eq_x_n+1_right}, $N-(c_n+d_n)^2>0$ and $c_{n+1}=c_n+d_n$ so $N>c_{n+1}^2$.

If $x_{n+1}$ is the left child of $x_n$ then, by Remark \ref{rem_e_c_d}, $N-c_{n+1}^2=\frac{N-c_n^2}{d_n}{d_{n+1}}>0$ and so $N>c_{n+1}^2$.

Proposition \ref{prop_N_c_n} is proved by induction. \hfill $\square$

\begin{coro} --- The LPP of $\alpha$ is an irrationality path and so $\alpha$ is irrational.
\label{coro_irrational}
\end{coro}

\emph{Proof}. --- For every non negative integer $n$, $\abs{c_n} < \sqrt{N}$ by Proposition \ref{prop_N_c_n} and $d_n \mid N-c_n^2$ by Proposition \ref{prop_e_c_d} so $0 < d_n \leqslant N-c_n^2 \leqslant N$. Hence, $(c_n)$ and $(d_n)$ are both bounded sequences of integers so there are two positive integers $n<m$ such that $c_n=c_m$ and $d_n=d_m$. Thus, $x_n=x_m$ and so the LPP of $\alpha$ is an irrationality path. \hfill $\square$

\section{Periodicity and symmetry}

Corollary \ref{coro_irrational} shows that the LPP of $\alpha$ contains at least two equal numbers $x_n$ et $x_m$. So, if we denote $n_0$ the index of the first number that appears at least twice in the LPP of $\alpha$ and $n_1$ the first index after $n_0$ such that $x_{n_0}=x_{n_1}$ then sequences $(c_n)_{n\geqslant n_0}$ and $(d_n)_{n\geqslant n_0}$ are periodic with least period $n_1-n_0$ and the LPP of $\alpha$ is periodic from $x_{n_0}$. We have seen in previous examples it seems that $n_0=0$. Thus, the LPP of $\alpha$ seems to be purely periodic, i.e., it seems there is a positive integer $T$ such that $x_{n+T}=x_n$ for every non negative integer $n$. We will now show it.

\begin{lem} --- Every number in the Calkin-Wilf tree of $\alpha$ has a unique representation in the form 
\[\alpha=\frac{\varepsilon\sqrt{N}+c}{d}\]
where $\varepsilon\in\{-1,1\}$, $c$ is a rational integer and $d$ is a positive integer.
\label{lem_uniqueness}
\end{lem}

\emph{Proof}. --- The existence of such a representation is guaranteed by Proposition \ref{prop_e_c_d}. Assume there are numbers $\varepsilon$ and $\varepsilon'$ in $\{-1,1\}$, rational integers $c$ and $c'$ and positive integers $d$ and $d'$ such that 
\[\dfrac{\varepsilon\sqrt{N}+c}{d}=\dfrac{\varepsilon'\sqrt{N}+c'}{d'}.\]
Thus, $(d'\varepsilon-d\varepsilon')\sqrt{N}=dc'-d'c$. But, by Corollary \ref{coro_irrational} (if one takes $p=0$ and $q=1$), $\sqrt{N}$ is irrational so $d'\varepsilon-d\varepsilon'=0$ and $dc'-d'c=0$. It follows $d=\abs{d\varepsilon'}=\abs{d'\varepsilon}=d'$, so $c=c'$ and $\varepsilon=\varepsilon'$ which proves the uniqueness. \hfill $\square$

\begin{lem} --- Let $x_{n}$ and $x_{m}$ be  two numbers in the LPP of $\alpha$. Assume that $x_{n+1}$ is the right child of $x_{n}$ and $x_{m+1}$ is the left child of $x_m$. Then, $x_{n+1} \neq x_{m+1}$.
\label{lem_only_one_parent}
\end{lem}

\emph{Proof}. --- Assume by contradiction $x_{n+1}=x_{m+1}$. By \eqref{eq_x_n+1_right}, \eqref{eq_x_n+1_left} and Lemma \ref{lem_uniqueness}, 
\[c_n+d_n=c_m-\frac{N-c_m^2}{d_m} \qquad \text{and} \qquad d_n=\frac{(c_m+d_m)^2-N}{d_m}=2c_m+d_m-\frac{N-c_m^2}{d_m}.\]
Thus, $c_n+2c_m+d_m=c_m$, i.e., $c_n=-(c_m+d_m)$. It follows that $c_n^2=(c_m+d_m)^2>N$ since $N-(c_m+d_m)^2<0$ by \eqref{eq_x_n+1_left}. This contradicts Proposition \ref{prop_N_c_n} so $x_{n+1} \neq x_{m+1}$. \hfill $\square$

\begin{prop} --- The LPP of $\alpha$ is purely periodic. 
\end{prop}

\emph{Proof}. --- Assume that the first index $n_0$ such that $x_{n_0}$ appears at least twice in the LPP of $\alpha$ is not $0$ and denote $n_1$ the second index in the path such that $x_{n_1}=x_{n_0}$. By Lemma \ref{lem_only_one_parent}, $x_{n_0}$ and $x_{n_1}$ are both either left children or right children of their parents. Furthermore, the maps $x\mapsto x+1$ and $x\mapsto \frac{x}{x+1}$ are clearly injective so $x_{n_0-1}=x_{n_1-1}$ which contradicts the minimality of $n_0$. Thus, $n_0=0$ and so the path is purely periodic. \hfill $\square$. 

\bigskip

In the following, we denote $T$ the least period of the LPP of $\alpha$, i.e., the least positive integer such that $x_{T}=\alpha$. Note that $T\geqslant 2$ since right and left children of an irrational number $x$ cannot be equal to $x$. One can observe in the above Calkin-Wilf tree of $\sqrt{5}$ some kind of symmetry in the pattern between $x_0$ and $x_{T}$. This is what we will now investigate. 

\begin{lem} --- Let $n$ be a positive integer. If $x_n$ is the right child of $x_{n-1}$ then 
\begin{equation}
(c_n-d_n)^2<N \qquad \text{and} \qquad x_{n-1}=\dfrac{\sqrt{N}+c_n-d_n}{d_n}.
\label{eq_left_parent}
\end{equation}
If $x_n$ is the left child of $x_{n-1}$ then 
\begin{equation}
(c_n-d_n)^2>N \qquad \text{and} \qquad x_{n-1}=\dfrac{\sqrt{N}+c_n+\frac{N-c_n^2}{d_n}}{\frac{(c_n-d_n)^2-N}{d_n}}.
\label{eq_right_parent}
\end{equation}
\label{lem_reverse}
\end{lem}

\emph{Proof}. --- If $x_n$ is the right child of $x_{n-1}$ then $x_n=x_{n-1}+1$ so 
\[x_{n-1}=x_n-1=\dfrac{\sqrt{N}+c_n-d_n}{d_n}.\]
Thus, by Lemma \ref{lem_uniqueness}, $c_{n-1}=c_n-d_n$ so, by Proposition \ref{prop_N_c_n}, $(c_n-d_n)^2<N$.

If $x_n$ is the left child of $x_{n-1}$ then $x_n=\frac{x_{n-1}}{x_{n-1}+1}$ so
\[x_{n-1}=\dfrac{x_n}{1-x_n}=\dfrac{\sqrt{N}+c_n}{d_n-\sqrt{N}-c_n}=\dfrac{(\sqrt{N}+c_n)(d_n+\sqrt{N}-c_n)}{(d_n-c_n)^2-N} = \dfrac{d_n\sqrt{N}+c_nd_n+N-c_n^2}{(c_n-d_n)^2-N}.\]

Since $d_n>0$, Lemma \ref{lem_uniqueness} yields $(c_n-d_n)^2-N>0$ and
\[x_{n-1}=\dfrac{\sqrt{N}+c_n+\frac{N-c_n^2}{d_n}}{\frac{(c_n-d_n)^2-N}{d_n}}\]
as announced. \hfill $\square$

\bigskip

Recall that if $\alpha$ and $\beta$ are two rational numbers, the algebraic conjugate of $x:=\alpha\sqrt{N}+\beta$ is  
\[x^*:=-\alpha\sqrt{N}+\beta.\]
Since $\sqrt{N}$ is irrational, rational numbers $\alpha$ and $\beta$ are unique as in Lemma \ref{lem_uniqueness} so $x^*$ is well-defined. Moreover, it is well-known that for every rational numbers $\gamma$ and $\delta$, 
\begin{equation}
(\gamma x + \delta)^*= \gamma x^*+\delta \qquad \text{ and } \qquad \left(\frac{1}{x}\right)^*=\frac{1}{x^*}.
\label{eq_x^*}
\end{equation}

\begin{rem} --- By Lemma \ref{lem_uniqueness}, if there are positive integers $m$ and $n$ such that $x_m=-x_n^*$ then $c_m=-c_n$ and $d_m=d_n$.
\label{rem_-x_n^*}
\end{rem}

\begin{prop} --- If $m\in\{1, 2, ..., T\}$ is such that $x_m=-\alpha^*$ then for every $n\in\{0, 1, ...,  m\}$, $x_{m-n}=-x_n^*$ and $x_{n+1}$ and $x_{m-n}$ are both either left children or right children of their parents.
\label{prop_symmetry}
\end{prop}

\emph{Proof}. --- Suppose that $m\in\{1, 2, ..., T\}$ is an index such that $x_m=-\alpha^*$.

Then, by definition, $x_m=-x_0^*$. Moreover, $(c_m-d_m)^2=(-c_0-d_0)^2=(c_0+d_0)^2$ thus $(c_0+d_0)^2$ and $(c_m-d_m)^2$ are both either greater or less than $N$ and so, by \eqref{eq_x_n+1_right}, \eqref{eq_x_n+1_left} and Lemma \ref{lem_reverse}, $x_{1}$ and $x_{m}$ are both either left children or right children of their parents.

Assume the property is true for a certain integer $n\in\{0, 1, ...,  m-1\}$. Then, $x_{m-n}=-x_{n}^*$ so $c_{m-n}=-c_n$ and $d_{m-n}=d_n$.

If $x_{n+1}$ is the right child of $x_n$ and $x_{m-n}$ is the right child of $x_{m-n-1}$ then, by \eqref{eq_left_parent} and \eqref{eq_x_n+1_right}, 
\[x_{m-n-1}=\dfrac{\sqrt{N}+c_{m-n}-d_{m-n}}{d_{m-n}}=\dfrac{\sqrt{N}-c_{n}-d_{n}}{d_{n}}=-(x_n+1)^*=-x_{n+1}^*.\]
If $x_{n+1}$ is the left child of $x_n$ and $x_{m-n}$ is the left child of $x_{m-n-1}$ then, by \eqref{eq_right_parent} and \eqref{eq_x_n+1_left},
\[x_{m-n-1}=\dfrac{\sqrt{N}+c_{m-n}+\frac{N-c_{m-n}^2}{d_{m-n}}}{\frac{(d_{m-n}-c_{m-n})^2-N}{d_{m-n}}}=\dfrac{\sqrt{N}-c_{n}+\frac{N-c_{n}^2}{d_{n}}}{\frac{(d_{n}+c_{n})^2-N}{d_{n}}}=\dfrac{\sqrt{N}-\left(c_{n}-\frac{N-c_{n}^2}{d_{n}}\right)}{\frac{(c_{n}+d_{n})^2-N}{d_{n}}}=-x_{n+1}^*.\]
Moreover, $(c_{m-n-1}-d_{m-n-1})^2=(-c_{n+1}-d_{n+1})^2=(c_{n+1}+d_{n+1})^2$ so, as previously, $x_{n+2}$ and $x_{m-n-1}$ are both either right children or left children of their parents. 

Thus, Proposition \ref{prop_symmetry} is proved by induction. \hfill $\square$

\begin{coro} --- Assume $x_m=-\alpha^*$ for a certain positive integer $m$. 
\begin{enumerate}
\item The list $(c_0, c_1, ..., c_{m-1},c_{m})$ is antisymmetric and the list $(d_0, d_1, ..., d_{m-1}, d_{m})$ is symmetric. In particular, if $m=2k$ is even then $c_{k}=0$ and $x_{k}=\frac{\sqrt{N}}{d_{k}}$.
\item The steps are symmetric that is to say the list $(s_0, s_1, ..., s_{m-1})$ is a palindrome: for every integer $j\in\{0, 1, ..., m-1\}$, $s_{m-1-j}=s_j$.
\item For every $n\in\{0, 1, ..., T\}$, there exists $k(n)\in\{0, 1, ... ,T\}$ such that $-x_n^*=x_{k(n)}$.
\end{enumerate}
\end{coro}

\emph{Proof}. --- Statements 1 and 2 immediately follow from Proposition \ref{prop_symmetry} and Lemma \ref{lem_uniqueness}. Proposition \ref{prop_symmetry} also yields statement 3 if $n\in\{0, 1, ..., m\}$. Moreover, if one considers the sequence $(x_n')$ related to the LPP of $x_m=-\alpha^*$, clearly $x_n'=x_{m+n}$ for every non negative $n$. Thus, since $x_{T-m}'=x_T=-x_m^*=-(x_0')^*$, by Proposition \ref{prop_symmetry}, for every $n\in\{0, 1, ..., T-m\}$, $-(x_{n}')^*=x_{T-m-n}'$, i.e., for every $n\in\{m, m+1, ..., T\}$, $-x_{n}^*=x_{T-(n-m)}$. \hfill $\square$

\begin{rem} \hspace{1cm}
\begin{enumerate}
\item There are numbers $\alpha$ such that $-\alpha^*$ does not appear in the LPP of $\alpha$. Such is the case, for instance, of $\frac{\sqrt{34}+1}{3}$ whose LPP is:
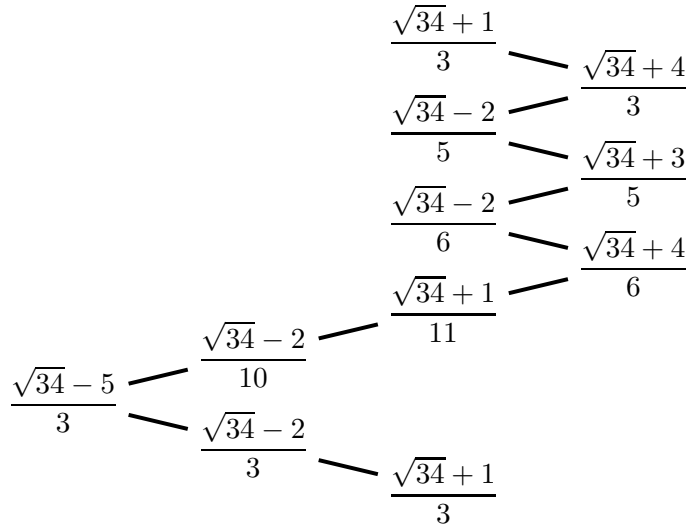
\begin{center}
\begin{tikzpicture}[xscale=1,yscale=1]
\node (n1) at ({0},{0}) {$\dfrac{\sqrt{34}+1}{3}$};
\node (n2) at ({2.5},{-0.6}) {$\dfrac{\sqrt{34}+4}{3}$};
\node (n3) at ({0},{-1.2}) {$\dfrac{\sqrt{34}-2}{5}$};
\node (n4) at ({2.5},{-1.8}) {$\dfrac{\sqrt{34}+3}{5}$};
\node (n5) at ({0},{-2.4}) {$\dfrac{\sqrt{34}-2}{6}$};
\node (n6) at ({2.5},{-3}) {$\dfrac{\sqrt{34}+4}{6}$};
\node (n7) at ({0},{-3.6}) {$\dfrac{\sqrt{34}+1}{11}$};
\node (n8) at ({-2.5},{-4.2}) {$\dfrac{\sqrt{34}-2}{10}$};
\node (n9) at ({-5},{-4.8}) {$\dfrac{\sqrt{34}-5}{3}$};
\node (n10) at ({-2.5},{-5.4}) {$\dfrac{\sqrt{34}-2}{3}$};
\node (n11) at ({0},{-6}) {$\dfrac{\sqrt{34}+1}{3}$};
\node (n12) at ({7},{-6}) {};
\draw[line width=1.5pt] (n1)--(n2);
\draw[line width=1.5pt] (n2)--(n3);
\draw[line width=1.5pt] (n3)--(n4);
\draw[line width=1.5pt] (n4)--(n5);
\draw[line width=1.5pt] (n5)--(n6);
\draw[line width=1.5pt] (n6)--(n7);
\draw[line width=1.5pt] (n7)--(n8);
\draw[line width=1.5pt] (n8)--(n9);
\draw[line width=1.5pt] (n9)--(n10);
\draw[line width=1.5pt] (n10)--(n11);
\end{tikzpicture}
\end{center}
\item If $p=0$, i.e., $\alpha=\frac{\sqrt{N}}{q}$ then $-\alpha^*=\alpha$ so $x_T=-\alpha^*$ and $m=T$ is the least index such that $x_m=-\alpha^*$. Moreover, in this case, the list $(s_0, s_1, ..., s_{T-1})$ is palindromic.
\end{enumerate}
\end{rem}

\section{Relationship with continued fractions}

The proof of periodicity of sequence $(x_n)$ has strong similarities to the proof of periodicity for the continued fraction expansion of quadratic surds due to Lagrange (see \cite{Lag70} and \cite[pp. 41-44]{RS92}). Furthermore, periodicity and symmetry of sequence $(s_n)$ remind properties of the continued fraction expansion of square roots established by Legendre (see \cite[pp. 42-47]{Leg08} and \cite[p. 47]{RS92}) and of reduced quadratic surds due to Galois (see \cite{Gal29} and \cite[p. 45]{RS92}). We will now explicit this similitude and show that the previous results lead to Legendre and Galois theorems.

\subsection{Introduction}

We define a subsequence $(y_n)$ of $(x_n)$ by putting $y_0=x_0=\alpha$ and, for every non negative integer $n$, $y_{n+1}$ is the first number after $y_n$ in the LPP of $\alpha$ where the direction changes, i.e., the first number $x_{\varphi(n)}$ after $y_n$ such that $s_{\varphi(n)-1} \neq s_{\varphi(n)}$. For every positive integer $n$, we say that $y_n$ is LR (left to right) if $s_{\varphi(n)-1}=\ell$ and $s_{\varphi(n)}=r$ and $y_n$ is RL (right to left) if $s_{\varphi(n)-1}=r$ and $s_{\varphi(n)}=\ell$. Furthermore, for every non negative integer $n$, we denote $t_n$ the number of steps between $y_n$ and $y_{n+1}$ in the LPP of $\alpha$, i.e., if $y_n=x_j$ and $y_{n+1}=x_{k}$ then $t_n=k-j$. Thus, for every positive integer $n$, if $y_n$ is RL then $t_n$ counts steps left and if $y_n$ is LR then $t_n$ counts steps right.

For example, if $\alpha=\frac{\sqrt{19}+4}{3}$ then one has the following simplified LLP:

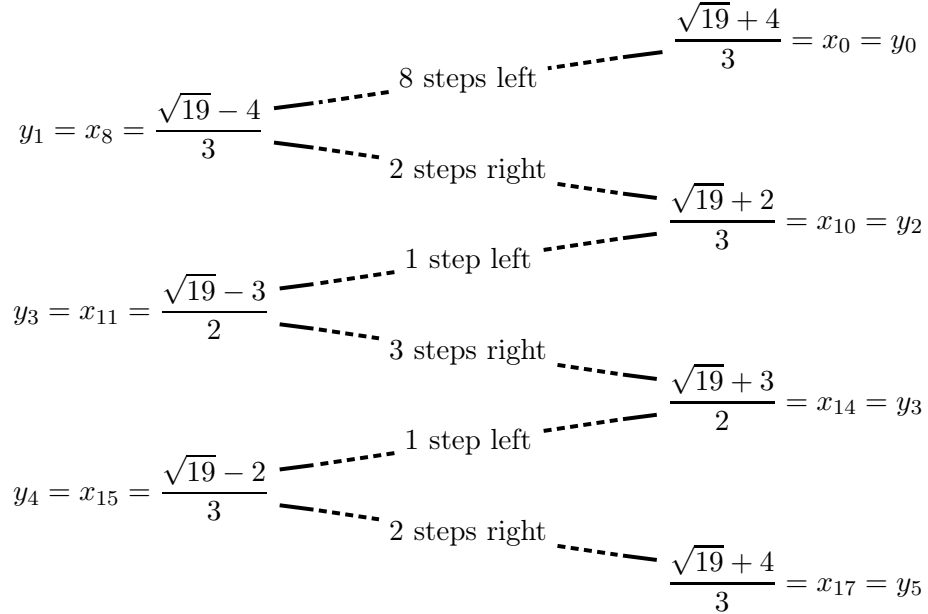
\begin{center}
\begin{tikzpicture}[xscale=1,yscale=1] 
\node (n1) at ({4.3},{0}) {$\dfrac{\sqrt{19}+4}{3}=x_0=y_0$};
\node[inner sep=0.5,outer sep=0.5] (n1i) at ({2},{-0.3}) {};
\node (n2) at ({0},{-0.6}) {$8$ steps left};
\node[inner sep=0.5,outer sep=0.5] (n2i) at ({-2},{-0.9}) {};
\node (n3) at ({-4.3},{-1.2}) {$y_1=x_8=\dfrac{\sqrt{19}-4}{3}$};
\node[inner sep=0.5,outer sep=0.5] (n3i) at ({-2},{-1.5}) {};
\node (n4) at ({0},{-1.8}) {$2$ steps right};
\node[inner sep=0.5,outer sep=0.5] (n4i) at ({2},{-2.1}) {};
\node (n5) at ({4.3},{-2.4}) {$\dfrac{\sqrt{19}+2}{3}=x_{10}=y_2$};
\node[inner sep=0.5,outer sep=0.5] (n5i) at ({2},{-2.7}) {};
\node (n6) at ({0},{-3}) {$1$ step left};
\node[inner sep=0.5,outer sep=0.5] (n6i) at ({-2},{-3.3}) {};
\node (n7) at ({-4.3},{-3.6}) {$y_3=x_{11}=\dfrac{\sqrt{19}-3}{2}$};
\node[inner sep=0.5,outer sep=0.5] (n7i) at ({-2},{-3.9}) {};
\node (n8) at ({0},{-4.2}) {$3$ steps right};
\node[inner sep=0.5,outer sep=0.5] (n8i) at ({2},{-4.5}) {};
\node (n9) at ({4.3},{-4.8}) {$\dfrac{\sqrt{19}+3}{2}=x_{14}=y_3$};
\node[inner sep=0.5,outer sep=0.5] (n9i) at ({2},{-5.1}) {};
\node (n10) at ({0},{-5.4}) {$1$ step left};
\node[inner sep=0.5,outer sep=0.5] (n10i) at ({-2},{-5.7}) {};
\node (n11) at ({-4.3},{-6}) {$y_4=x_{15}=\dfrac{\sqrt{19}-2}{3}$};
\node[inner sep=0.5,outer sep=0.5] (n11i) at ({-2},{-6.3}) {};
\node (n12) at ({0},{-6.6}) {$2$ steps right};
\node[inner sep=0.5,outer sep=0.5] (n12i) at ({2},{-6.9}) {};
\node (n13) at ({4.3},{-7.2}) {$\dfrac{\sqrt{19}+4}{3}=x_{17}=y_{5}$};
\draw[line width=1.5pt] (n1)--(n1i);
\draw[line width=1.5pt, densely dashed] (n1i)--(n2);
\draw[line width=1.5pt, densely dashed] (n2)--(n2i);
\draw[line width=1.5pt] (n2i)--(n3);
\draw[line width=1.5pt] (n3)--(n3i);
\draw[line width=1.5pt, densely dashed] (n3i)--(n4);
\draw[line width=1.5pt, densely dashed] (n4)--(n4i);
\draw[line width=1.5pt] (n4i)--(n5);
\draw[line width=1.5pt] (n5)--(n5i);
\draw[line width=1.5pt, densely dashed] (n5i)--(n6);
\draw[line width=1.5pt, densely dashed] (n6)--(n6i);
\draw[line width=1.5pt] (n6i)--(n7);
\draw[line width=1.5pt] (n7)--(n7i);
\draw[line width=1.5pt, densely dashed] (n7i)--(n8);
\draw[line width=1.5pt, densely dashed] (n8)--(n8i);
\draw[line width=1.5pt] (n8i)--(n9);
\draw[line width=1.5pt] (n9)--(n9i);
\draw[line width=1.5pt, densely dashed] (n9i)--(n10);
\draw[line width=1.5pt, densely dashed] (n10)--(n10i);
\draw[line width=1.5pt] (n10i)--(n11);
\draw[line width=1.5pt] (n11)--(n11i);
\draw[line width=1.5pt, densely dashed] (n11i)--(n12);
\draw[line width=1.5pt, densely dashed] (n12)--(n12i);
\draw[line width=1.5pt] (n12i)--(n13);
\end{tikzpicture}
\end{center}
The first few terms of $(y_n)$ are
\[\dfrac{\sqrt{19}+4}{3},~\dfrac{\sqrt{19}-4}{3},~\dfrac{\sqrt{19}+2}{3},~\dfrac{\sqrt{19}-3}{2}, ~\dfrac{\sqrt{19}+3}{2},~\dfrac{\sqrt{19}-2}{3},~\dfrac{\sqrt{19}+4}{3}, ...\]
and $y_1$ is LR, $y_2$ is RL, $y_3$ is LR and so on.
The first few terms of $(t_n)$ are $8,~2,~1,~3,~1,~2, ...$

\bigskip

Note that for every non negative integer $k$, $x_{k+1}=x_k+1$ if $x_{k+1}$ is the right child of $x_k$ and $\frac{1}{x_{k+1}}=\frac{1}{x_k}+1$ if $x_{k+1}$ is the left child of $x_k$. Thus, for every non negative integer $n$,
\begin{equation}
y_{n+1}=y_n+t_n \text{ if $y_{n+1}$ is RL} \qquad \text{and} \qquad \dfrac{1}{y_{n+1}}=\dfrac{1}{y_n}+t_n \text{ if $y_{n+1}$ is LR}.
\label{eq_expression_y_n_t_n}
\end{equation}

\begin{prop} --- For every positive integer $n$,  
\label{prop_LR}
\[0<-y_n^*<1<y_n \text{ if $y_n$ is RL} \qquad \text{and} \qquad 0 < y_n < 1 < -y_n^* \text{ if $y_n$ is LR.}\]
\end{prop}

\emph{Proof}. --- Let $n$ be a positive integer and write $y_n=\frac{\sqrt{N}+c}{d}$ as in Proposition \ref{prop_e_c_d}.

If $y_n$ is RL then, by \eqref{eq_left_parent} and \eqref{eq_x_n+1_left}, $(c-d)^2<N<(c+d)^2$ so $-2cd<2cd$ which implies $c>0$. Therefore, $d-c\leqslant \abs{d-c} < \sqrt{N}<\abs{c+d}=c+d$ so $y_n=\frac{\sqrt{N}+c}{d}>1$ and $-y_n^*=\frac{\sqrt{N}-c}{d}<1$. Moreover, by Proposition \ref{prop_N_c_n}, $c^2<N$ so $-y_n^*>0$.  

If $y_n$ is LR then, by \eqref{eq_x_n+1_right} and  \eqref{eq_right_parent}, $(c+d)^2<N<(c-d)^2$ so $2cd<-2cd$ which implies $c<0$. Therefore, $c+d \leqslant \abs{c+d} < \sqrt{N}<\abs{c-d}=d-c$ so, similarly, $-y_n^*>1$ and $0<y_n<1$. \hfill $\square$

\bigskip

If $x$ is a real number, we denote $\lfloor x \rfloor$ the integer part of $x$ that is to say the greatest integer not exceeding $x$.

\begin{coro} --- For every non negative integer $n$, 
\label{coro_t_n}
\[t_{n}=\lfloor -y_{n}^* \rfloor \text{ if $y_{n+1}$ is RL} \qquad \text{and} \qquad t_{n}=\left\lfloor -\frac{1}{y_{n}^*}\right\rfloor \text{ if $y_{n+1}$ is LR.}\]
\end{coro}

\emph{Proof}. --- Let $n$ be a non negative integer. 

If $y_{n+1}$ is RL then, by \eqref{eq_expression_y_n_t_n}, $y_{n+1}=y_n+t_n$ so, by \eqref{eq_x^*}, $-y_n^*=-y_{n+1}^*+t_n$. But, by Proposition \ref{prop_LR}, $0<-y_{n+1}^*<1$ so $\lfloor -y_n^* \rfloor = t_n$.

If $y_{n+1}$ is LR then, by \eqref{eq_expression_y_n_t_n}, $\frac{1}{y_{n+1}}=\frac{1}{y_n}+t_n$ so, by \eqref{eq_x^*}, $-\frac{1}{y_n^*}=-\frac{1}{y_{n+1}^*}+t_n$. But, by Proposition \ref{prop_LR}, $-y_{n+1}^*>1$ thus $0<-\frac{1}{y_{n+1}^*}<1$ and so $\left\lfloor -\frac{1}{y_n^*} \right\rfloor = t_n$. \hfill $\square$

\subsection{Square roots and Legendre theorem}

Let $[a_0, a_1, a_2, ...]$ be the (regular) continued fraction expansion of $\alpha$ that is to say $a_0$, $a_1$, $a_2$, ... are all positive integers (since $\alpha>0$) such that
\[\alpha=a_0+\dfrac{1}{a_1+\dfrac{1}{a_2+\dfrac{1}{\ddots}}}.\]
The sequence $(a_n)_{n\geqslant 0}$ is defined by the following recursion:
\[\begin{cases} \zeta_0=\alpha \text{ and } a_0=\lfloor  \zeta_0 \rfloor \\
  \text{for every integer } n\geqslant 0,~ \zeta_{n+1}=\dfrac{1}{ \zeta_n-a_n} \text{ and } a_{n+1}=\lfloor  \zeta_{n+1} \rfloor \end{cases}.\]
  The auxiliary sequence $( \zeta_n)_{n\geqslant 0}$ is the sequence of the complete quotients.

 \begin{prop} --- Assume $p=0$ and $N>q^2$. Then, for every non negative integer $n$, 
 \begin{enumerate}
 \item[$\bullet$] $\zeta_n=- y_n^*$ if $n$ is even and $\zeta_n=-\frac{1}{ y_n^*}$ if $n$ is odd;
 \item[$\bullet$] $a_n=t_n$.
 \end{enumerate}
 \label{prop_continued_fraction_sqrt}
 \end{prop} 

\emph{Proof}. --- By Remark \ref{rem_positive_path}, the LPP of $\alpha$ begins with a step right. Thus, by Corollary \ref{coro_t_n}, for every non negative integer, $t_{n}=\lfloor -y_{n}^* \rfloor$ if $n$ is even and $t_{n}=\left\lfloor -\frac{1}{y_{n}^*}\right\rfloor$ if $n$ is odd. 

Since $\zeta_0=\frac{\sqrt{N}}{q}=-\left(\frac{\sqrt{N}}{q}\right)^*=- y_0^*$ and $a_0=\left\lfloor \frac{\sqrt{N}}{q} \right\rfloor= t_0$, the property is true for $n=0$.

Assume it is true for a certain non negative integer $n$.

If $n$ is even then $y_{n+1}=y_n+t_n$ so, by \eqref{eq_x^*} and \eqref{eq_expression_y_n_t_n}, 
\[\zeta_{n+1}= \frac{1}{\zeta_n-a_n}=\frac{1}{-y_n^*-t_n}=-\frac{1}{(y_n+t_n)^*}=-\frac{1}{y_{n+1}^*}\]
and, since $n+1$ is odd,  $a_{n+1}=\lfloor  \zeta_{n+1} \rfloor = \left\lfloor -\frac{1}{y_{n+1}^*} \right\rfloor=t_{n+1}$.

If $n$ is odd then $y_{n+1}=\frac{y_n}{t_ny_n+1}$ so, by \eqref{eq_expression_y_n_t_n} and  \eqref{eq_x^*}, 
\[\frac{1}{\zeta_{n+1}}=\zeta_n-a_n=-\frac{1}{y_n^*}-t_n=-\left(\frac{1}{y_n}+t_n\right)^*=-\frac{1}{y_{n+1}^*}\]
so $\zeta_{n+1}=-y_{n+1}^*$ and, since $n+1$ is even,  $a_{n+1}=\lfloor  \zeta_{n+1} \rfloor=\lfloor -y_{n+1}^* \rfloor=t_{n+1}$.

Proposition \ref{prop_continued_fraction_sqrt} is proved by induction. \hfill $\square$

\bigskip

We deduce the following generalization of a theorem due to Legendre (see \cite[p. 47]{Leg08} for original statement and \cite[p. 47]{RS92} for generalization).

\begin{thm} --- Let $R>1$ be a rational number that is not the square of a rational number. The continued fraction expansion of $\sqrt{R}$ has the form
\[\sqrt{R}=\left[a_0, \overline{a_1, a_2, ..., a_2, a_1, 2a_0}\right].\]
\end{thm}

\emph{Proof}. --- Let us write $R=\frac{f}{g}$ in reduced form and put $N=fg$ and $q=g$. Then, $\sqrt{R}=\frac{\sqrt{N}}{q}$ and, since $R>1$, $N>q^2$. Moreover, $N$ cannot be a perfect square. Otherwise, since $f$ and $g$ are coprime, $f$ and $g$ would be perfect squares and so $R$ would be the square of a rational number.

Denote $y_m$ the last term of $(y_n)$ that appears before $x_{T}=\sqrt{R}$ in the LPP of $\sqrt{R}$. Since $(s_n)$ is $T$-periodic and $(s_0, s_1, ..., s_{T-1})$ is a palindrome, there are $t_0$ steps right between $y_m$ and $x_{T}$ and $t_0$ steps right between $x_{T}$ and $y_{m+1}$ so $t_{m}=2t_0$. It follows that the sequence $(t_n)$ is periodic from $t_1$ and is in the form
\begin{equation}
t_0,~t_1,~t_2, ~\cdots, t_2,~t_{1},~2t_0,~t_1,~t_2, , ~\cdots, ~t_2, ~t_{1}, 2t_0, ...
\end{equation}
One concludes by Proposition \ref{prop_continued_fraction_sqrt}. \hfill $\square$

\begin{rem} \hspace{1cm}
\begin{enumerate}
\item Since $a_0=\lfloor \sqrt{N} \rfloor$, it is clear that $\sqrt{N}+\lfloor \sqrt{N} \rfloor=\left[\overline{2a_0, a_1, a_2, ..., a_2, a_1}\right]$ so this expansion is purely periodic. It is a special case of a theorem due to Galois (see below).
\item It is well known (see \cite[p. 40]{RS92}) that for every non negative integer $n$, 
\[ \zeta_n=-\dfrac{q_{n-2}\sqrt{N}-p_{n-2}}{q_{n-1}\sqrt{N}-p_{n-1}}\]
where $\frac{p_n}{q_n}$ is the reduced form of the convergent $R_n=[a_0, a_1, ..., a_n]$ (with convention $p_{-2}=q_{-1}=0$ and $q_{-2}=p_{-1}=1$). It follows from Proposition \ref{prop_continued_fraction_sqrt} that, if $\alpha=\frac{\sqrt{N}}{q}$ with $N>q^2$ then, for every non negative integer $n$,
\[y_n=\dfrac{q_{n-2}\sqrt{N}+p_{n-2}}{q_{n-1}\sqrt{N}+p_{n-1}} \text{ if $n$ is even} \qquad \text{ and } \qquad y_n=\dfrac{q_{n-1}\sqrt{N}+p_{n-1}}{q_{n-2}\sqrt{N}+p_{n-2}} \text{ if $n$ is odd.}\]
\end{enumerate}
\end{rem}

\subsection{Reduced quadratic surds and Galois theorem}

Let us recall that an irrational number $x$ is said to be a reduced quadratic surd if $x$ is a root of a quadratic polynomial with integer coefficients $P$ such that $x>1$ and its conjugate $x^*$ satisfies $-1<x^*<0$. For example, the golden ratio $\varphi=\frac{\sqrt{5}+1}{2}$ and the number $\frac{\sqrt{19}+4}{3}$ studied above are reduced quadratic surds.

\begin{lem} --- Let $x$ be a reduced quadratic surd. Then, $x$ can be written in the form 
\[x=\dfrac{\sqrt{N}+p}{q}\]
where $N$ is positive integer which is not a perfect square, $p$ is a positive integer such that $p^2<N$ and $q$ is a positive integer which divides $N-p^2$. Moreover, $(p-q)^2<N<(p+q)^2$.
\label{lem_reduced_form}
\end{lem}

\emph{Proof}. --- Let $P=aX^2+bX+c$ be a quadratic polynomial with integer coefficients such that $P(x)=0$. Without loss of generality, one can assume that $a>0$. Put $\Delta:=b^2-4ac$ the discriminant of $P$ so that $x=\frac{-b\pm\sqrt{\Delta}}{2a}$. Since $x$ is an irrational number, $\Delta$ is not a perfect square and $x^*$ is the other root of $P$. Since $x>0$ and $x^*<0$, $\frac{c}{a}=xx^*$ is negative so $c<0$. Thus, $\Delta>b^2$ and so $-b-\sqrt{\Delta}<0$. It follows that $x=\frac{-b+\sqrt{\Delta}}{2a}$. Put $N=\Delta$, $p=-b$ and $q=2a$ so that $x=\frac{\sqrt{N}+p}{q}$, $x^*=\frac{-\sqrt{N}+p}{q}$ and $N>p^2$. Then, $N-p^2=-4ac=(-2c)q$ so $q$ divides $N-p^2$. Since $x^*>-1$, $-\sqrt{N}+p>-q$ thus $p+q>\sqrt{N}$ and, since $x>1$, $\sqrt{N}+p>q$ thus $q-p<\sqrt{N}<q+p$. This clearly implies $p>0$ and $(p+q)^2>N$. Moreover, if $p\leqslant q$ then $0\leqslant q-p < \sqrt{N}$ so $(p-q)^2<N$ and if $p>q$ then $0<p-q<p<\sqrt{N}$ so $(p-q)^2<N$. \hfill $\square$

\bigskip

From now on, we assume that $x$ is a reduced quadratic surd. Then, $\alpha:=-\frac{1}{x^*}$ is also a reduced quadratic surd so one can write $\alpha=\frac{\sqrt{N}+p}{q}$ as in Lemma \ref{lem_reduced_form}. We consider sequences $(\zeta_n)$ and $(a_n)$ related to the continued fraction expansion of $x$ and sequences $(x_n)$, $(y_n)$, $(s_n)$ and $(t_n)$ related to the LPP of $\alpha$ (and not those related to $x$). 

\begin{prop} ---  For every non negative integer $n$, 
 \begin{enumerate}
 \item[$\bullet$] $\zeta_n=- \frac{1}{y_n^*}$ if $n$ is even and $\zeta_n=-y_n^*$ if $n$ is odd;
 \item[$\bullet$] $a_n=t_{n}$.
 \end{enumerate}
 \label{prop_continued_fraction_reduced}
 \end{prop} 

\emph{Proof}. --- By Lemma \ref{lem_reduced_form}, $(p+q)^2 > N$ so, by \eqref{eq_x_n+1_left}, the LPP of $\alpha$ begins with a step left. Hence, by Corollary \ref{coro_t_n}, for every integer $n\geqslant 0$, $t_n=\left\lfloor -\frac{1}{y_n^*} \right\rfloor$ if $n$ is even and $t_n=\lfloor -y_n^* \rfloor$ is $n$ is odd.

Since, $y_0=\alpha=-\frac{1}{x^*}$, $-\frac{1}{y_0^*}=x=\zeta_0$ and $a_0=\lfloor \zeta_0 \rfloor=\left\lfloor -\frac{1}{y_0^*} \right\rfloor=t_0$.

Assume that the property is true for a certain integer $n\geqslant 0$.

If $n$ is even then by \eqref{eq_x^*} and \eqref{eq_expression_y_n_t_n},
\[\frac{1}{\zeta_{n+1}}=\zeta_n-a_n=-\frac{1}{y_n^*}-t_n=-\left(\frac{1}{y_{n+1}}+t_n\right)^*=-\frac{1}{y_{n+1}^*}\]
so $\zeta_{n+1}=-y_{n+1}^*$ and, since $n+1$ is odd, $a_{n+1}=\lfloor  \zeta_{n+1} \rfloor=\lfloor -y_{n+1}^* \rfloor=t_{n+1}$.

Similarly, if $n$ is odd, 
\[\zeta_{n+1}= \frac{1}{\zeta_n-a_n}=\frac{1}{-y_n^*-t_n}=-\frac{1}{y_{n+1}^*}\]
and $a_{n+1}=\lfloor  \zeta_{n+1} \rfloor = \left\lfloor -\frac{1}{y_{n+1}^*} \right\rfloor=t_{n+1}$.

So, Proposition \ref{prop_continued_fraction_reduced} is proved by induction. \hfill $\square$

\bigskip

We deduce the following theorem due to Galois (see \cite{Gal29} and \cite[p. 45-46]{RS92}).

\begin{thm} --- If $x$ is a reduced quadratic surd then its continued fraction expansion is purely periodic. Moreover, if this expansion is $x=\left[\overline{a_0, a_1, \ldots, a_{m-1}}\right]$
then $-\frac{1}{x^*}=\left[\overline{a_{m-1}, a_{m-2}, \ldots, a_{0}}\right]$.
\label{theo_Galois}
\end{thm}

\emph{Proof}. --- The sequence $(x_n)$ is periodic with period $T$ thus $x_T=\alpha=\frac{\sqrt{N}+p}{q}$. By Lemma \ref{lem_reduced_form}, $(p-q)^2 < N < (p+q)^2$ so, by \eqref{eq_left_parent} and \eqref{eq_x_n+1_left}, the path changes direction at $x_T$ from right to left. Thus, there is a positive integer $m$ such that $\alpha=x_{T}=y_m$ and $y_m$ is RL. The sequence $(t_n)$ is so periodic with period $m$ and one deduces from Proposition \ref{prop_continued_fraction_reduced} that $(a_n)$ is also periodic with period $m$. Moreover, since $y_m$ is RL, the pattern between $y_0$ and $y_m$ begins with a step left and ends with a step right so $m$ has to be even.

Denote $(y_n')$ and $(t_n')$ the sequences related to the LPP of $x$. Let us prove by induction that for every non negative $n$, $y_{n}'=-\frac{1}{y_{m-n}^*}$ and $t_n'=t_{m-n-1}$.

Since $y_m=\alpha=-\frac{1}{x^*}$, $-\frac{1}{y_m^*}=x=y_0'$. Moreover, since $x$ is a reduced quadratic surd, the LPP of $x$ also begins with a step left so $t_0'=\left\lfloor -\frac{1}{(y_0')^*} \right\rfloor = \lfloor y_m \rfloor$. Since $y_m$ is RL, $y_{m-1}$ is LR and so, by Proposition \ref{prop_LR}, $0<y_{m-1}<1$. Moreover, by \eqref{eq_expression_y_n_t_n}, $y_{m}=y_{m-1}+t_{m-1}$ so $\lfloor y_m \rfloor = t_{m-1}$ and thus $t_0'=t_{m-1}$.

Assume that the property is true for a certain non negative integer $n$.

If $n$ is even then $y_{n+1}'$ is LR so, by \eqref{eq_expression_y_n_t_n} and \eqref{eq_x^*},
\[\dfrac{1}{y_{n+1}'}=\dfrac{1}{y_n'}+ t_n'=-y_{m-n}^*+t_{m-1-n}=-(y_{m-n}-t_{m-n-1})^*=-y_{m-n-1}^*\]
because, since $m$ is even, $m-n$ is also even and thus $y_{m-n}=y_{m-n-1}+t_{m-n-1}$. Hence, $y_{n+1}'=-\frac{1}{y_{m-n-1}^*}$. Moreover, since $n+1$ and $m-n-1$ are both odd, by Corollary \ref{coro_t_n} and \eqref{eq_expression_y_n_t_n},
\[t_{n+1}'=\lfloor -(y_{n+1}')^*\rfloor = \left\lfloor \frac{1}{y_{m-n-1}} \right\rfloor=\left\lfloor \frac{1}{y_{m-n-2}}+t_{m-n-2} \right\rfloor = t_{m-n-2}\]
because $m-n-2$ is even so, by Proposition \ref{prop_LR}, $y_{m-n-2}>1$. 

If $n$ is odd then $y_{n+1}'$ is RL so, by \eqref{eq_expression_y_n_t_n} et \eqref{eq_x^*},
\[y_{n+1}'=y_n'+t_n'=-\frac{1}{y_{m-n}^*}+t_{m-n-1}=-\left(\dfrac{1}{y_{m-n}}-t_{m-n-1}\right)^*=-\dfrac{1}{y_{m-n-1}^*}\]
because $m-n-1$ is even. Moreover, since $n+1$ and $m-n-1$ are both even, by Corollary \ref{coro_t_n} and \eqref{eq_expression_y_n_t_n},
\[t_{n+1}'= \left\lfloor -\frac{1}{(y_{n+1}')^*} \right\rfloor = \left\lfloor y_{m-n-1} \right\rfloor = \left\lfloor y_{m-n-2}+t_{m-n-2} \right\rfloor  = t_{m-n-2}\]
because $m-n-2$ is odd so $0< y_{m-n-2}<1$ by Proposition \ref{prop_LR}. Thus, $t_{n+1}'=t_{m-(n+1)-2}$.

The property is proved by induction and the form of continued fraction expansion of $-\frac{1}{x^*}$ follows immediately from Proposition \ref{prop_continued_fraction_reduced}. \hfill $\square$

\begin{rem} \hspace{1cm}

\begin{enumerate}
\item The complete Galois theorem states the converse is true: if the continued fraction expansion of a number $x$ is purely periodic then $x$ is a reduced quadratic surd.
\item With previous notations, it follows from Proposition \ref{prop_continued_fraction_reduced} and Theorem \ref{theo_Galois} that sequences $(a_n)$ and $(t_n')$ related to $x$ have same period $m$ and, for every integer $n\in\{0, 1, ..., m-1\}$, $a_n=t_{m-n-1}'$. Hence, if one wants to read  the continued fraction expansion of $x$ directly from the LPP of $x$, it suffices to take values of $(t_n')$ in the reverse order. For example, the LPP drawn at the beginning of this section shows that $\frac{\sqrt{19}+4}{3}=\left[\overline{2, 1, 3, 1, 2, 8}\right]$.
\item The period $m$ yield by the LPP is not always the least one. For example, the LPP of $\frac{\sqrt{37}+5}{3}$ leads to $m=6$ and $\frac{\sqrt{37}+5}{3}=\left[ \overline{3,2,1,3,2,1}\right]$ so the least period is actually $3$.
\end{enumerate}
\end{rem}

\noindent \textbf{Acknowledgement}. --- The author would like to thank Daniel Duverney for his help with the exposition of this article. 

\bibliographystyle{amsalpha}
\bibliography{biblio}

\end{document}